\documentclass[10pt,a4paper,twoside]{amsart}

\usepackage{amsfonts, amssymb, amsmath, amsthm}
\usepackage{mathrsfs} 
\usepackage{latexsym}
\usepackage{enumerate}

\usepackage{url}

\usepackage[colorlinks=true]{hyperref}
\hypersetup{citecolor=blue, linkcolor=blue}

\newtheorem{thm}{Theorem}[section]
\newtheorem{lem}[thm]{Lemma}
\newtheorem{cor}[thm]{Corollary}

\def\Ocal{\mathcal{O}}
\def\1{1\!\!\!1}

\newcommand{\myK}{\mathbb{K}}

\title{Explicit count of integral ideals of an
  imaginary quadratic field%
}
\author{Olivier Ramar\'e}

\begin{document}
 \address[O. Ramar\'e]{CNRS/ Institut de Math\'ematiques de Marseille, Aix 
 Marseille Universit\'e, U.M.R. 7373, Site Sud, Campus de Luminy, Case 907, 
 13288 
 Marseille Cedex 9, France.}
 \email{olivier.ramare@univ-amu.fr}



\maketitle

\begin{abstract}
  We provide explicit bounds for the number of integral ideals of
  norms at most~$X$ is $\mathbb{Q}[\sqrt{d}]$ when $d <0$ is a
  fundamendal discriminant with an error term of size
  $\Ocal(X^{1/3})$. In particular, we prove that, when $\chi$
  is the non-principal character modulo~$3$ and $X\ge1$, we have
  $\sum_{n\le X}(\1\star\chi)(n)
    =
    \frac{\pi X}{3\sqrt{3}}
    +\Ocal^*(
    1.94\,X^{1/3})$, and that , when $\chi$
  is the non-principal character modulo~$4$ and $X\ge1$, we have
  $\sum_{n\le X}(\1\star\chi)(n)
    =
 \frac{\pi X}{4}
    +\Ocal^*(
    1.4\,X^{1/3}
    )$.
\end{abstract}


\section{Introduction and results}

\subsubsection*{General context}
Let $\myK$ be a number field, of degree $n_\myK$, discriminant
$\Delta(\myK)$, associated Dedekind zeta-function $\zeta_\myK$ of
residue $\kappa_\myK$ at one. Counting the number of integral ideals
of norm below some bound is a fundamental question that has been
addressed by numerous authors. The explicit angle has been less
popular and three pieces of works emerge: the paper \cite{Debaene*19}
by K.~Debaene,
the PhD memoir \cite{Sunley*73} by J.~Sunley and its upgraded version
\cite{Lee*22} by E.S.~Lee. The first goes by lattice point counting,
gets a dependence on the size of order $x^{1-1/n_\myK}$
and treats the dependence in the
field very finely. This approach is reused in
\cite{Gun-Ramare-Sivaraman*22b} to enumerate integral ideals in ray
classes.  The second and third approaches follow the analytic
treatment proposed by E.~Landau: they get a better dependence on the size of
order $x^{1-2/(n_\myK+1)}$ but only rely on the discriminant of the
field, an invariant which is notoriously large. The constants obtained
have the clear advantage of being explicit but they remain
gigantic: when $n_\myK=2$, the error term in \cite{Lee*22} reads
\begin{equation*}
  \Ocal^*\biggl(8.81\cdot 10^{11}|\Delta_\myK|^{\frac{1}{n_\myK+1}}
  (\log|\Delta_\myK|)^{n_\myK-1} x^{1-\frac{2}{n_\myK+1}}\biggr).
\end{equation*}
We aim here at being less demanding in generality but to gain in
numerical precision.
\subsubsection*{Our results for imaginary quadratic number fields}
Let $d$ be a squarefree integer. We associate to this integer
its fundamental discriminant defined by
\begin{equation}
  \label{defDelta}
  \Delta_d=
  \begin{cases}
    d&\text{when $d\equiv 1[4]$},\\
    4d&\text{when $d\equiv 2,3[4]$}.
  \end{cases}
\end{equation}
The associated character is given in terms of the Kronecker symbol by
the formula~$\chi(n)=(\frac{\Delta_d}{n})$.
\begin{thm}
  \label{Main}
  When $X\ge \max(|\Delta_d|, 2c_0(d))$ and $d$ is a negative squarefree integer, we have
  \begin{equation*}
    \sum_{n\le X}(\1\star\chi)(n)
    =
    XL(1,\chi)
    +\frac{1}{2|\Delta_d|}
       \sum_{1\le r\le |\Delta_d|}r\chi(r)
    +\Ocal^*\bigl(
    0.76\,L(1,\chi)c_0(d)X^{1/3}
    \bigr)
  \end{equation*}
  where
  \begin{equation}
    \label{defc0d}
    c_0(d)=\max(c(3/4), c(5/4))^{2/3},
  \end{equation}
  and
  \begin{equation}
    \label{eq:1}
    c(3/4)=\max_{M\ge1}\sum_{m\le M}
    \frac{(\1\star\chi)(m)}{m^{3/4}M^{1/4}L(1,\chi)},
    \
    c(5/4)=\frac{M^{1/4}}{L(1,\chi)}\sum_{m\ge M}
      \frac{(\1\star\chi)(m)}{m^{5/4}}.
    \end{equation}
    When $X\ge \max(130^2|\Delta_d|,2c_0(d))$, the constant $0.76$ may be replaced by $0.67$.
\end{thm}
Lemmas~\ref{DedekindSum} and~\ref{DedekindSum2} propose upper bounds
for $c(3/4)$ and $c(5/4)$. It is in particular proved that
$\min(c(3/4), c(5/4))\ge 4$. As the reader sees, imposing larger bounds reduces only
marginally the final constant. It may still be of interest for very
small values of $d$, where the range in $X$ can be completed by direct computations.

\subsubsection*{Special cases}
We start by some numerical verification done with a GP-Pari script.
\begin{thm}
  \label{Verif4}
    Let $\chi$ be the non principal character modulo~4. For $X\in[1,
    10^8]$, we have
    \begin{equation*}
    \sum_{n\le X}(\1\star\chi)(n)
    =
    \frac{\pi X}{4}
    +\Ocal^*\bigl(
    2.08\,X^{1/4}
    \bigr).
  \end{equation*}
\end{thm}
The constant in the big-O seems to increase slightly when $X$ increases.
\begin{thm}
  \label{Verif3}
    Let $\chi$ be the non principal character modulo~3. For $X\in[1,
    10^8]$, we have
    \begin{equation*}
    \sum_{n\le X}(\1\star\chi)(n)
    =
    \frac{\pi X}{3\sqrt{3}}
    +\Ocal^*\bigl(
    1.63\,X^{1/4}
    \bigr).
  \end{equation*}
\end{thm}

\begin{cor}
  Let $\chi$ be the non principal character modulo~4. For $X\ge 1$, we have
  \begin{equation*}
    \sum_{n\le X}(\1\star\chi)(n)
    =
    \frac{\pi X}{4}
    +\Ocal^*\bigl(
    1.4\,X^{1/3}
    \bigr).
  \end{equation*}
\end{cor}
\begin{cor}
  Let $\chi$ be the non principal character modulo~3. For $X\ge 1$, we have
  \begin{equation*}
    \sum_{n\le X}(\1\star\chi)(n)
    =
    \frac{\pi X}{3\sqrt{3}}
    +\Ocal^*\bigl(
    1.94\,X^{1/3}
    \bigr).
  \end{equation*}
\end{cor}
\begin{cor}
  Let $\chi$ be the quadratic character of
  $\mathbb{Q}[\sqrt{d}]$, where~$-19\le d \le -1$.
  For $X\ge 68$, we have
  \begin{equation*}
    \sum_{n\le X}(\1\star\chi)(n)
    =
     X L(1,\chi)
    +\frac{1}{2|\Delta|}
       \sum_{1\le r\le |\Delta_d|}r\chi(r)
    +\Ocal^*\bigl(
    3.4\,X^{1/3}
    \bigr).
  \end{equation*}
\end{cor}

\subsubsection*{Methodology}

E.~Landau's approach in \cite{Landau*17h} (See also \cite[Satz 210]{Landau*49}) relies on several ingredients, but the first and
main one is the functional equation of the associated Dedekind zeta
function. E.~Landau sends the line of integration to $\Re s=-1/2$,
then uses the functional equation to study the last integral.
This Landau's approach is described in modern language in~\cite{Lee*22}. 
The
process used has become known as the Vorono\"{\i} Summation
Formula(s), based on
\cite{Voronoi*04-1, Voronoi*04-2}, though this latter is more
commonly used for the divisor function.
In fact, though the papers of Vorono\"{\i} largely predates the ones
of Landau, and it cannot be assumed that Landau did not know of them,
Landau does not mention the Vorono\"{\i} approach, a surprising fact
as this author has most of the time been very prompt in explaining
the genesis of ideas. This absence may be due to
the combination of two facts:
Landau worked in great generality, with intricate
Gamma-factors, and a general view of the Vorono\"{\i}
process was missing at the time. 
This process is now well understood and is for instance well-documented in Chapter~10 of the book~\cite{Cohen*07} by H.~Cohen.
The addition of Vorono\"{\i} is to recognize the involved
Mellin transform as a Bessel function and to consider a functional transform of the initial
weight function, see Lemma~\ref{VoronoiSpe} below. We follow this approach here.

Two more ingredients are being used: a non-negative smoothing device
and an apriori trivial upper bound for the number of integral ideals
below some bound to avoid divisor functions of $x$ in the remainder
term, see Lemma~\ref{DedekindSum} and~\ref{DedekindSum2} below.

We do not examine what happens for the small values of the size with
respect to the discriminant.

\subsubsection*{Thanks}
Thanks are due to Olivier Bordell\`es for some helpful discussions on
this subject.

\section{On the Bessel functions}

\begin{lem}
  \label{derivBesselJ}
  When $a>0$, we have
  \begin{equation*}
    \int_0^X J_0(a\sqrt{t})dt =
    \frac{2\sqrt{X}}{a}J_1(a\sqrt{X})
  \end{equation*}
  and
  \begin{equation*}
    \int_0^X tJ_0(a\sqrt{t})dt =
    \frac{4X}{a^2}J_2(a\sqrt{X})
    +
     \frac{2X^{3/2}}{a}J_1(a\sqrt{X}).
  \end{equation*}
\end{lem}

\begin{proof}
  Indeed, when $\nu\ge1$, we have $J_\nu(t)'=
  J_{\nu-1}(t)-J_\nu(t)\nu/t$.
  Hence 
  \begin{equation*}
    (\sqrt{t}J_1(a\sqrt{t}))'
    =\frac{1}{2\sqrt{t}}J_1(a\sqrt{t})
      +\frac{a}{2}\biggl(J_0(a\sqrt{t})-\frac{1}{a\sqrt{t}}J_1(a\sqrt{t})\biggr)
      =\frac{a}{2}J_0(a\sqrt{t})
  \end{equation*}
  which proves the first formula. For the second one, we notice
  similarly that
  \begin{align*}
    (t^{3/2}J_1(a\sqrt{t}))'
    &=\frac{3\sqrt{t}}{2}J_1(a\sqrt{t})
      +\frac{at}{2}\biggl(J_0(a\sqrt{t})-\frac{1}{a\sqrt{t}}J_1(a\sqrt{t})\biggr)
    \\& =\frac{a}{2}tJ_0(a\sqrt{t})
    -\sqrt{t}J_1(a\sqrt{t})
  \end{align*}
  and
  \begin{equation*}
    (tJ_2(a\sqrt{t}))'
    =
      J_2(a\sqrt{t})
      +\frac{a\sqrt{t}}{2}
      \biggl(
      J_1(a\sqrt{t})-\frac{2}{a\sqrt{t}}J_2(a\sqrt{t})
      \biggr)
      =\frac{a\sqrt{t}}{2}J_1(a\sqrt{t})
  \end{equation*}
  and therefore
  \begin{equation*}
    \frac{d}{dx}\biggl(
    \frac{at^{3/2}}{2}J_1(a\sqrt{t})
    +tJ_2(a\sqrt{t})
    \biggr)
    =\frac{a^2}{4}tJ_0(a\sqrt{t})
  \end{equation*}
\end{proof}

\begin{lem}
  \label{BoundsBessel}
  When $\nu>0$ and $x\ge 0$, we have
  \begin{equation*}
    \biggl|J_\nu(x)-\sqrt{\frac{2}{\pi x}}\cos(x-(2\nu+1)\tfrac{\pi}{4})\biggr|
    \le \frac{4|\nu^2-1/4|}{5x^{3/2}}.
  \end{equation*}
  When $\nu>1/2$ and $x\ge0$, we have
  $\displaystyle
    |x^2-\nu^2+\tfrac14|^{1/4}|J_\nu(x)|\le \sqrt{2/\pi}$.
\end{lem}

\begin{proof}
  The first inequality is given in a consequence of \cite[Theorem
  4]{Krasikov*14} by Krasikov while the
  second one comes from \cite[Theorem 3]{Krasikov*14}.
\end{proof}

\begin{lem}
  \label{BoundT}
  With 
  \begin{equation}
    \label{defT}
    T(z; a)=
    \frac{(1+z)J_2(a\sqrt{1+z})
      -J_2(a)}{z}
  \end{equation}
  we have, when $a\ge 4\pi$ and $z\in(0,1/3]$:
  \begin{equation*}
    |T(z,a)|\le
    \min\biggl(0.53\sqrt{a}, \frac{7}{3z\sqrt{a}}\biggr).
  \end{equation*}
  When $a\ge 130\cdot 4\pi$ and $z\in(0,1/10)$:
  \begin{equation*}
    |T(z,a)|\le
    \min\biggl(0.4\sqrt{a}, \frac{2.1}{z\sqrt{a}}\biggr).
  \end{equation*}
\end{lem}

\begin{proof}
  By the mean value theorem and when $z>0$, we find that
  \begin{align*}
    |T(z;a)|
    &\le
    \frac{a(\sqrt{1+z}-1)}{z}\max_{a\le t\le a\sqrt{1+z}}|J_2'(t)|
    + |J_2(a\sqrt{1+z})|
    \\&\le
    \frac{a}2\max_{a\le t\le a\sqrt{1+z}}
    \biggl(|J_1(t)|+\frac{2}{t}|J_2(t)|\biggr)
    + |J_2(a\sqrt{1+z})|.
  \end{align*}
  The map $x\mapsto |x^2-3/4|^{1/4}$ is non-decreasing when $x\ge
  \sqrt{3}/2$ and the map $x\mapsto |x^2-7/4|^{1/4}$ is also non-decreasing when $x\ge
  \sqrt{7}/2$. On assuming that $a\ge \sqrt{7}/{2}$,
  Lemma~\ref{BoundsBessel} thus gives us that
  \begin{equation*}
    \sqrt{\pi/2}|T(z;a)|\le
    \frac{a}{2|a^2-3/4|^{1/4}}
    + \frac{2}{|a^2-7/4|^{1/4}}.
  \end{equation*}
  When $a\ge 4\pi$, a rapid plot shows that
  $\sqrt{\pi/2}|T(z;a)|\le 0.66\sqrt{a}$, while, when $a\ge 130\cdot
  4\pi$, we find that $\sqrt{\pi/2}|T(z;a)|\le 0.5013\sqrt{a}$.
  This establishes the first bound. When $a$ is large, it is better to
  simply use
  \begin{align*}
    |T(z;a)|
    &\le\frac{|J_2(a\sqrt{1+z})|+
      |J_2(a)|}{z}+|J_2(a\sqrt{1+z})|
    \\&\le\frac{2}{z|a^2-7/4|^{1/4}}+\frac{1}{|a^2(1+z)-7/4|^{1/4}}
    \le \frac{7}{3z\sqrt{a}}
  \end{align*}
  which we prove first for $z\le 1/4$ by using
  $|a^2(1+z)-7/4|^{1/4}\ge |a^2-7/4|^{1/4}$ and then on discretizing
  the interval (such precision is not required, it only leads to be
  better looking estimate).
  When $z\le 1/10$ and $a\ge 130\cdot4\pi$, we find that $|T(z;a)|\le 2.1/(z\sqrt{a})$.
\end{proof}

\section{Some a priori estimates}

Let us start with a well-known estimate.
\begin{lem}
  \label{ZetaSum}
  When $s>0$ and $s\neq1$, we have
  \begin{equation*}
    \sum_{m\le M}\frac{1}{m^s}
    =\frac{M^{1-s}}{1-s}+\zeta(s)+\Ocal^*(1/M^s).
  \end{equation*}
\end{lem}

\begin{proof}
  Indeed, we find that
  \begin{align*}
    \sum_{m\le M}\frac{1}{m^s}
    &=s\int_1^M [t]\frac{dt}{t^{s+1}}+\frac{[M]}{M^s}
      =\frac{s}{s-1}-s\int_1^M\{t\}\frac{dt}{t^{s+1}}+\frac{M^{1-s}}{1-s}
      -\frac{\{M\}}{M^{s}}
    \\&=
    \frac{s}{s-1}-s\int_1^\infty\{t\}\frac{dt}{t^{s+1}}+\frac{M^{1-s}}{1-s}
    +s\int_M^\infty\{t\}\frac{dt}{t^{s+1}}
    -\frac{\{M\}}{M^{s}}
    \\&=
    \zeta(s)+\frac{M^{1-s}}{1-s}
    +s\int_M^\infty\{t\}\frac{dt}{t^{s+1}}
    -\frac{\{M\}}{M^{s}}.
  \end{align*}
  We finally check that
  \begin{equation*}
    s\int_M^\infty\{t\}\frac{dt}{t^{s+1}}
    -\frac{\{M\}}{M^{s}}
    =s\int_M^\infty(\{t\}-\{M\})\frac{dt}{t^{s+1}}
  \end{equation*}
  and the lemma follows readily.
\end{proof}
The next lemma gives an upper bound as well as a mean to approximate~$c(s)$, when $s\in(0,1)$.
\begin{lem}
  \label{DedekindSum}
  When $s\in(0,1)$, we have
  \begin{equation*}
    \sum_{n\le N}\frac{(\1\star\chi)(n)}{n^sN^{1-s}}
    \le
    \frac{L(1,\chi)}{1-s}
    +
    \frac{|\zeta(s)L(s,\chi)|}{N^{1-s}}
    +
    \frac{(\frac{1}{4}+|\zeta(s)|+\frac{5}{1-s})\Omega(\chi)}{\sqrt{N}}
  \end{equation*}
  where
  \begin{equation}
    \label{defOmegachi}
    \Omega(\chi)=\max_{L\ge1}\biggr|\sum_{\ell\le L}\chi(\ell)\biggr|.
  \end{equation}
  Moreover, the quantity considered is asymptotic to $L(1,\chi)/(1-s)$.
\end{lem}

\begin{proof}
  By using Lemma~\ref{ZetaSum} and the Dirichlet hyperbola formula, we find that the sum $S$ to be computed equals (with parameters $L\ge1$
  and $M\ge1$ such that $LM=N$)
  \begin{align*}
    S
    &=
      \sum_{\ell\le L}\frac{\chi(\ell)}{\ell^s}
      \biggl(\zeta(s)+\frac{(N/\ell)^{1-s}}{1-s}
      +\Ocal^*\biggl(\frac{\ell^s}{4N^s}\biggr)\biggr)
      +\sum_{m\le M}\frac{1}{m^s}
      \Ocal^*\bigl(2\Omega(\chi)/L^s\bigr)
    \\&=
    \zeta(s)\sum_{\ell\le L}\frac{\chi(\ell)}{\ell^s}
    +\frac{N^{1-s}}{1-s}\sum_{\ell\le L}\frac{\chi(\ell)}{\ell}
    +\Ocal^*\biggl(
    \frac{L}{4N^s}
    +2\Omega(\chi)\frac{1-s+M^{1-s}}{(1-s)L^s}
    \biggr)
    \\&=
    \zeta(s)L(s,\chi)
    +\frac{N^{1-s}L(1,\chi)}{1-s}
    \\&\qquad+\Ocal^*\biggl(
    \frac{N^{1-s}\Omega(\chi)}{(1-s)L}
    +\frac{L}{4N^s}
    +2\Omega(\chi)\frac{1-s+M^{1-s}}{(1-s)L^s}
    +|\zeta(s)|\frac{\Omega(\chi)}{L^s}
    \biggr)
  \end{align*}
  so that
  \begin{align*}
    S/N^{1-s}
    \le
    &\frac{|\zeta(s)L(s,\chi)|}{N^{1-s}}
    +\frac{L(1,\chi)}{1-s}
    +
    \frac{\Omega(\chi)}{(1-s)L}
    +\frac{L}{4N}
   \\& +2\Omega(\chi)\frac{1-s+M^{1-s}}{(1-s)L^sN^{1-s}}
    +|\zeta(s)|\frac{\Omega(\chi)}{L^sN^{1-s}}.
  \end{align*}
  The simplistic choice $L=M=N^{1/2}$ leads to
   \begin{align*}
    S/N^{1-s}
    \le
    \frac{|\zeta(s)L(s,\chi)|}{N^{1-s}}
    +\frac{L(1,\chi)}{1-s}
    +\frac{5\Omega(\chi)+\frac14}{(1-s)\sqrt{N}}
    +|\zeta(s)|\frac{\Omega(\chi)}{N^{1-s/2}}.
  \end{align*}
\end{proof}

Our final tool in this section is the next lemma which gives an upper
bound as well as a mean to approximate~$c(s)$, when $s>1$.
\begin{lem}
  \label{DedekindSum2}
  When $s>1$, we have
  \begin{equation*}
    N^{s-1}\sum_{n\ge N}\frac{(\1\star\chi)(n)}{n^s}
    \le
    \frac{L(1,\chi)}{s-1}
    +\frac{(3\zeta(s)+\frac{1}{s-1})\Omega(\chi)+\frac14}{\sqrt{N}}
     .
  \end{equation*}
  Moreover, the quantity considered is asymptotic to $L(1,\chi)/(s-1)$.
\end{lem}

\begin{proof}
  We proceed as in Lemma~\ref{DedekindSum}.
  On setting $S=\sum_{n\le N}(\1\star\chi)(n)/n^s$, we find that (with parameters $L\ge1$
  and $M\ge1$ such that $LM=N$)
  \begin{align*}
    S
    &=
      \sum_{\ell\le L}\frac{\chi(\ell)}{\ell^s}
      \biggl(\zeta(s)-\frac{(\ell/N)^{s-1}}{s-1}
      +\Ocal^*\biggl(\frac{\ell^s}{4N^s}\biggr)\biggr)
      +\sum_{m\le M}\frac{1}{m^s}
      \Ocal^*\bigl(2\Omega(\chi)/L^s\bigr)
    \\&=
    \zeta(s)\sum_{\ell\le L}\frac{\chi(\ell)}{\ell^s}
    -\frac{1}{(s-1)N^{s-1}}\sum_{\ell\le L}\frac{\chi(\ell)}{\ell}
    +\Ocal^*\biggl(
    \frac{L}{4N^s}
    +2\Omega(\chi)\frac{\zeta(s)}{L^s}
    \biggr)
    \\&=
    \zeta(s)L(s,\chi)
    -\frac{L(1,\chi)}{(s-1)N^{s-1}}
    +\Ocal^*\biggl(
    \frac{\Omega(\chi)}{(s-1)LN^{1-s}}
    +\frac{L}{4N^s}
    +3\Omega(\chi)\frac{\zeta(s)}{L^s}
    \biggr)
  \end{align*}
  so that
  \begin{align*}
    (\zeta(s)L(s,\chi)-S)N^{s-1}
    \le
    \frac{L(1,\chi)}{s-1}
    +
    \frac{\Omega(\chi)}{(s-1)L}
    +\frac{L}{4N}
    +3\Omega(\chi)\frac{\zeta(s)N^{s-1}}{L^s}
    .
  \end{align*}
  We take $L=M=N^{1/2}$ and get
   \begin{align*}
    (\zeta(s)L(s,\chi)-S)N^{s-1}
    \le
    \frac{L(1,\chi)}{s-1}
    +\frac{(3\zeta(s)+\frac{1}{s-1})\Omega(\chi)+\frac14}{\sqrt{N}}
     .
   \end{align*}
   (explain the shift $N+1\mapsto N$).
\end{proof}

\begin{proof}
  We readily find that
  \begin{equation*}
    \sum_{\ell\le L}\frac{\chi(\ell)}{\ell^{s}}
    =
    s\int_1^L\sum_{\ell\le t}\chi(\ell)\frac{dt}{t^{s+1}}
    +\frac{\sum_{\ell\le L}\chi(\ell)}{L^s}
  \end{equation*}
  and thus
  \begin{equation*}
    L^s\sum_{\ell\ge L}\frac{\chi(\ell)}{\ell^{s}}
    = sL^s\int_L^\infty\sum_{\ell\le t}\chi(\ell)\frac{dt}{t^{s+1}}
    -\sum_{\ell\le L}\chi(\ell)
  \end{equation*}
\end{proof}

\section{Around the Vorono\"{\i} Summation Formula}
Here is the Vorono\"{\i} Summation Formula we want to use.
\begin{lem}
  \label{VoronoiSpe}
  Let $f:[0,\infty)\mapsto\mathbb{C}$ be a function that we assume
  to have a finite number of simple discontinuities (i.e. with finite left
    and right limits), to be is piecewise $C^\infty$ and piecewise
    monotonic and such that
  $f$ and all its derivatives tend to zero faster than
    any power of~$x$ at infinity.
  We have
  \begin{align*}
    \sum_{n\ge1}(\1\star\chi)(n)f(n)
    =&L(1,\chi)\check{f}(0)
       +\frac{f(0)}{2|\Delta_d|}
       \sum_{1\le r\le |\Delta_d|}r\chi(r)
    \\&+
    \frac{2\pi}{\sqrt{|\Delta_d|}}
    \sum_{m\ge1}
    (\1\star\chi)(m)\int_0^\infty f(t)J_0\bigl(
    4\pi\sqrt{nt/|\Delta_d|}\bigr)dt.
  \end{align*}
\end{lem}

\begin{proof}
Let us denote by $\zeta_d$ the Dedekind zeta-function of the imaginary
quadratic field $\mathbb{Q}[\sqrt{d}]$. It satisfies the functional equation
\begin{equation}
  \label{eq:2}
  \gamma_d(1-s)\zeta_d(1-s)=  \gamma_d(s)\zeta_d(s)
  \quad\text{where}\quad
  \gamma_d(s)=\biggl(\frac{\sqrt{|\Delta_d|}}{2\pi}\biggr)^{s}\Gamma(s).
\end{equation}
We use \cite[Theorem 10.2.17]{Cohen*07} by H.~Cohen. The kernel to
  be considered is
  \begin{align*}
    K_d(x)
    &=\frac{1}{2i\pi}\int_{\frac12-i\infty}^{\frac12+i\infty}
    \biggl(\frac{\sqrt{|\Delta_d|}}{2\pi}\biggr)^{s}
    \biggl(\frac{\sqrt{|\Delta_d|}}{2\pi}\biggr)^{s-1}
    \frac{\Gamma(s)}{\Gamma(1-s)}x^{-s}ds
    \\
    &=
      \frac{1}{2i\pi} \frac{2\pi}{\sqrt{|\Delta_d|}}
        \int_{\frac12-i\infty}^{\frac12+i\infty}
    \frac{4^s\Gamma(s)}{\Gamma(1-s)}\biggl(\frac{16\pi^2x}{|\Delta_d|}\biggr)^{-s}ds.  
  \end{align*}
  The formula
  \begin{equation*}
    4^s\frac{\Gamma(s)}{\Gamma(1-s)}
    =\int_0^\infty t^{s-1}J_0(\sqrt{t})dt
  \end{equation*}
  gives us
  \begin{equation*}
    K_d(x)=
    \frac{2\pi}{\sqrt{|\Delta_d|}}
    J_0\bigl(
    4\pi\sqrt{x/|\Delta_d|}
    \bigr).
  \end{equation*}
  Concerning the value at 0, we use $\zeta(0)=-1/2$ and
  \begin{equation}
    \label{eq:3}
    L(0,\chi)=\frac{-1}{|\Delta_d|}
       \sum_{1\le r\le |\Delta_d|}r\chi(r)
     \end{equation}
     as per \cite[Corollary 10.3.2]{Cohen*07}.
\end{proof}

\begin{lem}
  \label{VoronoiSmooth}
  We have
  \begin{align*}
    \sum_{n\le X}\biggl(1-\frac{n}{X}\biggr)(\1\star\chi)(n)
    =&\frac{XL(1,\chi)}{2}
    +\frac{1}{2|\Delta_d|}
       \sum_{1\le r\le |\Delta_d|}r\chi(r)
    \\&+
    \frac{\sqrt{|\Delta_d|}}{2\pi}
    \sum_{m\ge1}
    \frac{(\1\star\chi)(m)}{m}J_2\bigl(
    4\pi\sqrt{mX/|\Delta_d|}\bigr).
  \end{align*}
\end{lem}

\begin{proof}
  By using Lemma~\ref{VoronoiSpe}, we readily find that
  \begin{align*}
    \sum_{n\le X}\biggl(1-\frac{n}{X}\biggr)(\1\star\chi)(n)
    =&\frac{XL(1,\chi)}{2}
    +\frac{1}{2|\Delta_d|}
       \sum_{1\le r\le |\Delta_d|}r\chi(r)
    \\&+
    \frac{2\pi}{\sqrt{|\Delta_d|}}
    \sum_{m\ge1}
    (\1\star\chi)(m)\int_0^X \biggl(1-\frac{t}{X}\biggr)J_0\bigl(
    4\pi\sqrt{mt/|\Delta_d|}\bigr)dt.
  \end{align*}
  By using Lemma~\ref{derivBesselJ}, we find that
  \begin{align*}
    \sum_{n\le X}\biggl(1-\frac{n}{X}\biggr)(\1\star\chi)(n)
    =&\frac{XL(1,\chi)}{2}
    +\frac{1}{2|\Delta_d|}
       \sum_{1\le r\le |\Delta_d|}r\chi(r)
    \\&+
    \frac{2\pi}{\sqrt{|\Delta_d|}}
    \sum_{m\ge1}
    (\1\star\chi)(m)\frac{|\Delta_d|}{4\pi^2 m}J_2\bigl(
    4\pi\sqrt{mX/|\Delta_d|}\bigr).
  \end{align*}
  Our lemma follows swiftly from this last expression.
\end{proof}

\section{Main engine}
\begin{lem}
  \label{FirstApprox}
  When $Y\in[0,X/3]$ and $X\ge |\Delta_d|$, we have
  \begin{multline*}
    \sum_{n\le X}(\1\star\chi)(n)
    +\sum_{X< n\le X+Y}\frac{X+Y-n}{Y}(\1\star\chi)(n)
    \\=
    \frac{(2X+Y)L(1,\chi)}{2}
    +\frac{1}{2|\Delta_d|}
       \sum_{1\le r\le |\Delta_d|}r\chi(r)
    +\Ocal^*\Bigl(
    0.36\,C_0(d)\sqrt{X|\Delta_d|/Y}
    \Bigr)
  \end{multline*}
  where $C_0(d)=L(1,\chi)c_0(d)$.
  When $Y\in[0,X/10]$ and $X\ge 130^2|\Delta_d|$, the constant $0.36$
  may be reduced to~$0.292$.
\end{lem}

\begin{proof}
  By using Lemma~\ref{VoronoiSmooth} twice, we find that
  \begin{align*}
    (1/Y)&\sum_{n\le X}
    \biggl[(X+Y)\biggl(1-\frac{n}{X+Y}\biggr)-
    X\biggl(1-\frac{n}{X}\biggr)\biggr](\1\star\chi)(n)
    \\&=\frac{(2X+Y)L(1,\chi)}{2}
    +\frac{1}{2|\Delta_d|}
       \sum_{1\le r\le |\Delta_d|}r\chi(r)
    \\&+
    \frac{\sqrt{|\Delta_d|}}{2\pi}
    \sum_{m\ge1}
    \frac{(\1\star\chi)(m)}{m}
    T(Y/X; 4\pi\sqrt{mX/|\Delta_d|}).
  \end{align*}
  We now majorize the last sum by appealing to Lemma~\ref{BoundT}. 
  Recall that we assume that $X\ge |\Delta_d|$ (resp. $x\ge
  130^2|\Delta_d|)$.
  We use the first
  estimate of Lemma~\ref{BoundT} when
  \begin{equation*}
    4\pi\sqrt{mX/|\Delta_d|}\le \frac{7X}{3Y\cdot 0.53}
    \quad
    \biggl(\text{resp. }
    4\pi\sqrt{mX/|\Delta_d|}\le \frac{2.1X}{Y\cdot 0.4}
    \biggr)
  \end{equation*}
  i.e. when
  \begin{equation*}
    m\ge M=\frac{49|\Delta_d|X}{144(0.53\pi)^2 Y^2}
    \quad
    \biggl(\text{resp. }
    m\ge M=\frac{2.1^2|\Delta_d|X}{16(0.4\pi)^2 Y^2}
    \biggr).
  \end{equation*}
  We thus get
  \begin{align*}
    \frac{\sqrt{|\Delta_d|}}{2\pi}
    &\sum_{m\ge1}
    \frac{(\1\star\chi)(m)}{m}
    |T(Y/X; 4\pi\sqrt{mX/|\Delta_d|})|
    \\&\le
    0.53\frac{|\Delta_d|^{1/4}X^{1/4}}{\sqrt{\pi}}
    \sum_{m\le M}
    \frac{(\1\star\chi)(m)}{m^{3/4}}
    +
    \frac{7X^{3/4}}{3Y}\frac{|\Delta_d|^{3/4}}{4\pi^{3/2}}
    \sum_{m > M}
    \frac{(\1\star\chi)(m)}{m^{5/4}}.
  \end{align*}
  We recall taht $C_0(d)=L(1,\chi)c_0(d)$ where $c_0(d)$ is being
  defined in~\eqref{defc0d}. On appealing to 
   Lemmas~\ref{DedekindSum} and~\ref{DedekindSum2}, this leads to
  \begin{align*}
    \frac{\sqrt{|\Delta_d|}}{2\pi}
    &\sum_{m\ge1}
    \frac{(\1\star\chi)(m)}{m}
    |T(Y/X; 4\pi\sqrt{mX/|\Delta_d|})|
    \\&\le
    0.53C_0(d)\frac{|\Delta_d|^{1/4}X^{1/4}}{\sqrt{\pi}}M^{1/4}
    +
    C_0(d)\frac{7X^{3/4}}{3Y}\frac{|\Delta_d|^{3/4}}{4\pi^{3/2}M^{1/4}}
    \\&\le
    C_0(d)\frac{X^{1/2}\sqrt{0.53}\sqrt{7/3}}{Y^{1/2}\pi}|\Delta_d|^{1/2}
    \le 0.36\,C_0(d)\sqrt{X|\Delta_d|/Y}.
  \end{align*}
  When $Y\in[0,X/10]$ and $X\ge 130^2|\Delta_d|$, the constant $0.36$
  can be replaced by an upper bound for $\sqrt{2.1\cdot 0.4}/\pi$, e.g.
  $0.292$. The proof is complete.
\end{proof}

\begin{proof}[Final proof]
  We use Lemma~\ref{FirstApprox} with $X-Y+Y$ and $X+Y$. Thus
  \begin{align*}
    \frac{-Y}{2}L(1,\chi)-0.36\,C_0(d)\sqrt{(X-Y)|\Delta_d|/Y}
    &\le S-XL(1,\chi)
      +\frac{1}{2|\Delta_d|}
       \sum_{1\le r\le |\Delta_d|}r\chi(r)
    \\&\le \frac{Y}{2}L(1,\chi)+0.36\,C_0(d)\sqrt{X|\Delta_d|/Y}.
  \end{align*}
  We select
  \begin{equation}
    \label{defY}
    Y = \biggl(\frac{0.36\,C_0(d)\sqrt{X}}{L(1,\chi)}\biggr)^{2/3}
  \end{equation}
  getting the error term
  \begin{equation*}
    X^{1/3}C_0(d)^{2/3}L(1,\chi)^{1/3}\biggr(\frac{0.36^{2/3}}{2}+0.36^{2/3}\biggr).
  \end{equation*}
  We find that $Y\le X/3$ when
  \begin{equation*}
    3^{3/2}\frac{0.36\,C_0(d)}{L(1,\chi)}\le X.
  \end{equation*}
  The theorem follows readily in that case. The case $X\ge
  130^2|\Delta_d|$ is treated similarly.
\end{proof}

\section{Computing $C_0(d)$}

We use the script
  \texttt{ConvolutionAndVoronoi-01.gp}
and the function $\texttt{run(1000000, d)}$ therein to build the next table.

\bigskip
\begin{tabular}{|c||c|c|c|}
  \hline
  $d$&$\Delta_d$&$\Omega(\chi)$&$L(1,\chi)c_0(d)\le$\cr
  \hline
  $-1$&$-4$&1&$2.04$\cr
  $-2$&$-8$&2&$2.89$\cr
  $-3$&$-3$&1&$1.58$\cr
  $-5$&$-20$&4&$3.66$\cr
  $-6$&$-24$&4&$3.35$\cr
  $-7$&$-7$&2&$3.09$\cr
  $-10$&$-40$&4&$2.60$\cr
  $-11$&$-11$&3&$2.48$\cr
  $-13$&$-52$&5&$2.30$\cr
  $-14$&$-56$&8&$4.40$\cr
  $-15$&$-15$&3&$4.21$\cr
  $-17$&$-68$&8&$4.01$\cr
  $-19$&$-19$&3&$1.90$\cr
  \hline
\end{tabular}


\end{document}